\documentclass[12pt]{amsart}
\usepackage{amsmath,amssymb}
\newtheorem{theorem}{Theorem}
\newtheorem{proposition}[theorem]{Proposition}

\newtheorem{example}[theorem]{Example}
\def\ds{\displaystyle}
\def\C{\Bbb C}
\def\D{\Bbb D}
\def\O{\mathcal O}
\def\R{\hskip0.5mm\text{Re}\hskip0.5mm}
\def\I{\hskip0.5mm\text{Im}\hskip0.5mm}
\def\eps{\varepsilon}

\title[\tiny Boundary behavior of invariant functions on planar domains]
{Boundary behavior of invariant functions on planar domains}

\author{Nikolai Nikolov}
\author{Maria Trybu\l{}a}
\author{Lyubomir Andreev}

\address{N. Nikolov: Institute of Mathematics and Informatics\\Bulgarian Academy
of Sciences\\ Acad. G. Bonchev 8, 1113 Sofia, Bulgaria\newline
\indent Faculty of Information Sciences\\
State University of Library Studies and Information Technologies\\
Shipchenski prohod 69A, 1574 Sofia,
Bulgaria}\email{nik@math.bas.bg}
\address{M. Trybu\l{}a: Faculty of Mathematics and Informatics\\Adam Mickiewicz University
\\Umultowska 87, 61-001 Pozna\'{n}, Poland}
\email{maria.h.trybula@gmail.com}
\address{L. Andreev: Institute of Mathematics and Informatics\\Bulgarian Academy of Sciences
\\Acad. G. Bonchev 8, 1113 Sofia, Bulgaria,} \email{lyubomir.andreev@math.bas.bg}



\thanks{M. Trybula and L. Andreev are partially supported by the
Bulgarian National Science Found under contract DFNI-I 02/14.}

\begin{document}

\begin{abstract} Precise behavior of the Carath\'eodory, Kobayashi and Bergman metrics and
distances near smooth boundary points of planar domains is found under different assumptions
of regularity.
\end{abstract}

\maketitle

The most precise known estimates on the boundary behavior of the Carath\'eodory, Kobayashi
and Bergman metrics ($\gamma_D,$ $\kappa_D$ and $\beta_D$) of a bounded strictly pseudoconvex domain
$D$ in $\C^n$ with $\mathcal C^3$- or $\mathcal C^4$-smooth boundary can be found in \cite{DF,Fu,Ma}
(for the $C^{2,\eps}$-smooth case see \cite{N2}).

Set $d_D(z)=\mbox{dist}(z,\partial D).$

For $D\subset\C$ let $\gamma_D(z)=\gamma_D(z;1),$ etc., and
$$m_D\in A=\{\gamma_D,\kappa_D,\beta_D/\sqrt 2,\sqrt{\pi K_D}\},$$
where $K_D$ is the Bergman kernel on the diagonal ($\beta_D$ is defined if $\C\setminus D$ is not polar).
Note that the functions in $A$ coincide if $D\subsetneq\C$ is simply connected.

In the plane the mentioned above results can be read as follows:
\vskip1mm

1. If $D\Subset\C$ has $C^{2,\eps}$-smooth boundary ($\eps\in(0,1]$), then
the function $\ds\frac{2m_Dd_D-1}{d_D^{\eps/2}}$ is bounded;

2. If $D\Subset\C$ has $C^4$-smooth boundary, then $\ds m_D-\frac{1}{2d_D}$ is
bounded.
\smallskip

Our first goal is to obtain the precise behavior of $m_D$ near
$\mathcal C^2$-smooth boundary points of domains in $\C.$

\begin{proposition}\label{cur} If $a$ is a $\mathcal C^2$-smooth boundary point
of a domain $D\subset\C,$ then
$$\lim_{z\to a}\left(m_D(z)-\frac{1}{2d_D(z)}\right)=\frac{\chi_D(a)}{4},$$
where $\chi_{\partial D}(a)$ is the signed curvature of $\partial D$ at $a.$
\end{proposition}

We have a weaker result in the $\mathcal C^1$-smooth case.

\begin{proposition}\label{jn} (see \cite{JN}) If $a$ is a $\mathcal C^1$-smooth
boundary point of a domain $D\subset\C,$ then
$$\lim_{z\to a}2m_D(z)d_D(z)=1.$$
\end{proposition}

The following intermediate result holds in the $\mathcal C^{1,\eps}$-smooth case.

\begin{proposition}\label{eps} If $\eps\in(0,1]$ and $a$ is a $\mathcal C^{1,\eps}$-smooth
boundary point of a domain $D\subset\C,$  then $\ds\frac{2m_Dd_D-1}{d_D^\eps}$ is a bounded
function near $a.$
\end{proposition}

This result is sharp as the next example shows.

\begin{example}\label{exa} (a) If $\eps\in(0,1),$ then $\ds f(z)=z-\frac{z^{1+\eps}}{4}$ is a
univalent function on $\D_1=\{z\in\C:|z-1|<1\},$ $D=f(\D_1)$ has $\mathcal C^{1,\eps}$-smooth
boundary and
$$\ds\lim_{u\to 0+}\frac{2m_D(u)d_D(u)-1}{d_D^\eps(u)}=\frac{\eps}{4}.$$

\noindent(b) The domain $D=\D\cup\{z\in\C:|x|<1,y>0\}$ is $\mathcal C^{1,1}$-smooth and
$$\liminf_{z\to 1}\left(2m_D(z)-\frac{1}{d_D(z)}\right)=0<\frac14
=\limsup_{z\to 1}\left(2m_D(z)-\frac{1}{d_D(z)}\right).$$
\end{example}

Combining the proofs of Propositions \ref{cur} and \ref{eps}, one may obtain the following result
in the $\mathcal C^{2,\eps}$-smooth case whose proof we omit.

\begin{proposition}\label{e2} If $\eps\in(0,1)$ and $a$ is a $\mathcal C^{2,\eps}$-smooth
boundary point of a domain $D\subset\C,$ then $\ds\frac{4m_D-2/d_D-\chi_D}{d_D^\eps}$
is a bounded function near $a,$ where $\chi_D(z)$ is the signed curvature of $\partial D$
at the closest point to $z.$
\end{proposition}

A similar example to Example \ref{exa} shows that this result is sharp.
\smallskip

Recall now that the boundary behavior of Kobayashi distance $k_D$ of a bounded strictly
pseudoconvex domain $D$ in $\C^n$ is known up to $\pm$ constant (see \cite{BB}). The same behavior have
the Carath\'eodory distance $c_D,$ the Lempert function $l_D$ and, in the $\mathcal C^{2,\eps}$-smooth
case, the Bergman distance $b_D$ divided by $\sqrt{n+1}$ (see \cite{N2}).

Sharper estimates of these invariants in the case of Dini-smooth bounded domains in $\C$
can be found in \cite{N1,NT} ($k_D=l_D$ if $D\subset\C$). Our next aim is to find
the precise behavior of $p_D\in\{c_D,k_D,b_D/\sqrt2\}$ near Dini-smooth boundary points
of domains in $\C$ (these distances coincide if $D\subsetneq C$ is simply connected).
Recall that Dini-smoothness is bet\-ween $\mathcal C^1$- and $\mathcal C^{1,\eps}$-smoothness.

Set $\ds s_D(z,w):=\sinh^{-1}\frac{|z-w|}{2\sqrt{d_D(z)d_D(w)}}$
\smallskip

\hskip2.5cm$\ds=\log\frac{|z-w|+\sqrt{|z-w|^2+4d_D(z)d(w)}}{2\sqrt{d_D(z)d_D(w)}}.$
\smallskip

Note that $s_D=p_D$ if $D\subset\C$ is a half-plane.

\begin{proposition}\label{dif} If $a$ is a Dini-smooth boundary point
of a domain $D\subset\C,$ then
$$\lim_{z,w\to a}(p_D(z,w)-s_D(z,w))=0.$$
\end{proposition}

When $s_D\to 0,$ Proposition \ref{dif} can be improved even in the
$\mathcal C^1$-smooth case.

\begin{proposition}\label{rat} (a) If $a$ is a $\mathcal C^1$-smooth
boundary point of a domain $D\subset\C,$ then
$$\lim_{\substack{z,w\to a\\z\neq w}}\frac{p_D(z,w)}{s_D(z,w)}=1.$$

\noindent (b) If $D\subset\C$ is a $\mathcal C^1$-smooth bounded domain, then
$$\lim_{z\to\partial D}\frac{p_D(z,w)}{s_D(z,w)}=1
\quad\mbox{ uniformly in }w\neq z.$$
\end{proposition}

\section{Definitions}

\noindent{\bf 1.} A boundary point $p$ of a domain $D\subset C$ is said
to be Dini-smooth if $\partial D$ near $p$ is given by a Dini-smooth
curve $\gamma:[0,1]\rightarrow\mathbb{C}$ with $\gamma'\neq 0$ (i.e.
$\ds\int_0^1\frac{\omega(t)}{t}dt<\infty,$ where $\omega$ is the
modulus of continuity of $\gamma'$). A domain in $\C$ is called
Dini-smooth if all its boundary points are Dini-smooth.
\smallskip

\noindent{\bf 2.} Let $D$ be a domain in $\C^n.$ The Kobayashi (pseudo)distance
$k_D$ is the largest (pseudo)distance not exceeding the Lempert function
$$l_D(z,w)=\inf\{\tanh^{-1}|\alpha|:\exists\varphi\in\O(\D,D)
\hbox{ with }\varphi(0)=z,\varphi(\alpha)=w\},$$ where $\D$ is the
unit disc. Note that $k_D$ is  the integrated form of the
Kobayashi (pseudo)metric
$$\kappa_D(z;X)=\inf\{|\alpha|:\exists\varphi\in\O(\D,D) \hbox{
with }\varphi(0)=z,\alpha\varphi'(0)=X\}.$$

Denote by $c_D$ the Carath\'eodory (pseudo)distance of $D:$
$$c_D(z,w)=\sup\{\tanh^{-1}|\psi(w)|:\psi\in\O(D,\D)\hbox{ with }\psi(z)=0\}.$$
Then $c_D\le k_D\le l_D.$

The Bergman distance $b_D$ of $D$ is the integrated form of the Bergman
metric $\beta_D,$ i.e.
$$b_D(z,w)=\inf_\gamma\int_0^1\beta_D(\gamma(t);\gamma'(t))dt,
\quad z,w\in D,$$
where the infimum is taken over all smooth curves $\gamma:[0,1]\to
D$ with $\gamma(0)=z$ and $\gamma(1)=w.$

Recall that
$$\beta_D(z;X)=\frac{M_D(z;X)}{\sqrt{K_D(z)}},\quad z\in D,\
X\in\C^n,$$ where
$$M_D(z;X) =\sup\{|f'(z)X|:f\in L^2_h(D),\ ||f||_{L^2(D)}\le 1,\ f(z)=0\}$$
and
$$K_D(z)=\sup\{|f(z)|^2:f\in L^2_h(D),\ ||f||_{L^2(D)}\le 1\}$$
is the Bergman kernel on the diagonal (we assume that $K_D(z)>0$).

Note that $M_D\le M_G$ and $K_D\le K_G$ if $G$ is a subdomain of $D$ but
$\beta_D$ does not share this property, in general.

It is well-known that $c_D\le b_D.$

We refer to \cite{JP} for other properties of the above invariants.

\section{Proofs}

\noindent{\it Proof of Proposition \ref{cur}.} For $z\in D$ near $a$ there is a unique
$\mathcal C^2$-smooth point $z'\in\partial D$ such that $|z'-z|=d_D(z).$ Denote by $n_{z'}$
the inner unit normal vector to $\partial D$ at $z'$ and set $\phi_z(\tau)=z'+\tau n_{z'}$
($\tau\in\C$), $D_z=\phi_z^{-1}(D).$

For any $\eps>0$ there exists $\delta>0$ such that if $|z-a|<\delta,$
then $$D^{\eps,\delta}:=D^\eps\cap\delta\D\subset D_z\subset
D^{-\eps}\cup(\C\setminus\delta\overline\D)=:D^{-\eps,\delta},$$
where $D^{\pm\eps}=\{\tau\in\C:2\R\tau>\chi^{\pm\eps}|\tau|^2\}$
and $\chi^{\pm\eps}=\chi_{\partial D}(a)\pm\eps.$
Note that $D^{\pm\eps}$ is a disc, the complement of a disc or a half-plane.

Since $\ds\beta_D=\frac{M_D}{\sqrt{K_D}},$ to prove the result, we may replace
$\beta_D/\sqrt 2$ in $A$ by $\sqrt{cM_D},$ where $c=\sqrt{\pi/2}.$ Then
$$m_{D^{-\eps,\delta}}(d_D(z))\le m_D(z)\le m_{D^{\eps,\delta}}(d_D(z)).$$

\noindent{\scshape Lemma L.} ({\it localization}) Let $\Pi$ be the upper half-plane. Then
$$0\le m_{\Pi\cap\D}(z)-m_{\Pi}(z)\le\frac{|z|}{1-|z|^2}.$$

\noindent{\it Subproof.} The map $\ds f:z\mapsto\left(\frac{z+1}{z-1}\right)^2$ transforms
conformally $\Pi\cap\D$ onto $\Pi.$ Hence
\begin{align*}m_{\Pi\cap\D}(z)-m_{\Pi}(z)&=\frac{|f'(z)|}{2\I f(z)}-\frac{1}{2\I z}
=\frac{1}{2\I z}\left(\frac{|1-z^2|}{1-|z|^2}-1\right)\\
&\le\frac{|z^2-|z|^2|}{2\I z(1-|z|^2)}=\frac{|z|}{1-|z|^2}.
\end{align*}

Lemma L easily implies that $m_{D^{\pm\eps,\delta}}(\tau)-m_{D^{\pm\eps}}(\tau)\to 0$
as $\tau\to 0.$ Simple calculations show that
$\ds m_{D^{\pm\eps}}(\tau)-\frac{1}{2\tau}\to\frac{\chi^{\pm\eps}}{4}$
as $\tau\to 0+.$ Thus
$$\frac{\chi^{-\eps}}{4}\le\liminf_{z\to a}\left(m_D(z)-\frac{1}{2d_D(z)}\right)\le
\limsup_{z\to a}\left(m_D(z)-\frac{1}{2d_D(z)}\right)\le\frac{\chi^\eps}{4}.$$
It remains to let $\eps\to 0.$
\smallskip

\noindent{\it Proof of Proposition \ref{eps}.} We may replace $D^{\pm\eps}$ from above by
$$D^{\pm\eps}=\{z\in\C:x>\pm c|y|^{1+\eps}\},$$
where $c>0$ is a constant.

If $\eps=1,$ it follows as before that
$$-\frac{c}{2}\le\liminf_{z\to a}\left(m_D(z)-\frac{1}{2d_D(z)}\right)\le
\limsup_{z\to a}\left(m_D(z)-\frac{1}{2d_D(z)}\right)\le\frac{c}{2}.$$

Otherwise, we may choose a conformal map $\eta_\pm$ from $D^{\pm\eps}$ onto
the right half-plane $\Pi$ that extends $\mathcal C^{1+\eps}$-smoothly at $0$
such that $\eta_\pm(0)=0$ and $\eta_\pm(\Bbb R^+)=\Bbb R^+.$
Using Lemma L, it is enough to show that
$$g(x)=\frac{2x m_{D^{\pm\eps}}(x)-1}{x^\eps}$$
is a bounded function on $(0,1)$ which follows from the equalities
$$g(x)=\frac{x\eta'_\pm(x)-\eta_\pm(x)}{x^\eps\eta_\pm(x)}
=\frac{x}{\eta_\pm(x)}\cdot\frac{\eta_\pm'(x)-\eta_\pm'(\xi_x)}{x^\eps},
\quad\xi_x\in(0,x).$$

\noindent{\it Proof of Example \ref{exa}.} (a) The univalence of $f$ follows from
the inequality $\R f'>0$ on $\D_1$ (cf. \cite[Proposition 1.10]{Pom}). The
$\mathcal C^{1,\eps}$-smoothness of $\partial D$ is clear.

Choose now $c>0$ such that $D\supset G=\{u+iv\in c\Bbb D:c^{1-\eps}u>|v|^{1+\eps}\}.$
Set $\ds\delta=\frac{1+\eps}{1-\eps}.$ It is not difficult to compute that
$$\lim_{u\to 0+}\frac{u-d_G(u)}{u^\delta}=
\frac{2}{c^2(1+\eps)^\delta(1-\eps)}$$
and hence
$$\lim_{u\to 0+}\frac{u-d_D(u)}{u^{1+\eps}}=0.$$

Let $x=f^{-1}(u).$ Since $m_D(u)f'(x)=m_{\D_1}(x),$ then
$$\lim_{u\to 0+}\hskip-0.2mm\frac{2m_D(u)d_D(u)-1}{d_D^\eps(u)}\hskip-0.2mm
=\hskip-0.2mm\lim_{x\to 0+}\hskip-0.2mm\frac{\ds\frac{2}{x(2-x)}\cdot\frac{4x-x^{1+\eps}+o(x^{1+\eps})}
{4-(1+\eps)x^\eps}-1}{x^\eps}\hskip-0.2mm=\hskip-0.2mm\frac{\eps}{4}.$$

\noindent (b) Let $u$ and $v$ be the respective
$\liminf$ and $\limsup.$ Then $u\ge 0$ by convexity and $v\le 1/4,$ since $D$ satisfies
the interior 1-ball condition. On the other hand, Lemma L easily implies that if
$a,b\in\partial D\setminus\{\pm 1\}$ are such that $|\R a|=1$ and $|b|=1,$ then $\ds \lim_{z\to a}m_D(z)=0$
and $\ds\lim_{z\to b}m_D(z)=1/4$ which shows that $u\le 0$ and $v\ge 1/4.$
\smallskip

\noindent{\it Proof of Proposition \ref{dif} for $c_D$ and $k_D.$}
We shall use arguments from the proof of \cite[Proposition 5]{N1}.

We may find a Dini-smooth Jordan curve $\zeta$ such that
$\zeta=\partial D$ near $a$ and $D\subset\zeta_{\mbox{ext}}.$
Let $b\not\in\overline{\zeta_{\mbox{ext}}},$
$\ds\varphi(z)=\frac{1}{z-b}$ and $G=\zeta_{\mbox{ext}}\cup\{0\}.$
Let $\psi:G\to\D$ be a Riemann map. It extends to a $C^1$-diffeomorphism from
$\overline G$ to $\overline \D$ (cf. \cite[Theorem
3.5]{Pom}). Let $\theta$ maps conformally $\D$ onto a half-plane $\Pi.$
Setting $\eta=\theta\circ\psi\circ\varphi,$ then
$$c_D(z,w)\ge c_\Pi(\eta(z),\eta(w))=s_\Pi(\eta(z),\eta(w)).$$

Since $\ds\lim_{z\to a}\frac{d_\Pi(\eta(z))}{d_D(z)}=|\eta'(a)|,$
it follows that
$$\lim_{z,w\to a}(s_D(z,w)-s_\Pi(\eta(z),\eta(w))=0.$$
Hence
$$\liminf_{z,w\to a}(c_D(z,w)-s_D(z,w))\ge 0.$$

It remains to show that
$$\limsup_{z,w\to a}(k_D(z,w)-s_D(z,w))\le 0.$$
For this, we choose a Dini-smooth simply connected domain
$F\subset D$ such that $\partial F=\partial D$ near $a.$
Then we may proceed similarly to above.
\smallskip

In the next two proofs we shall use the quasi-hyperbolic distance $h_D,$ i.e.
the integrated form of $1/d_D.$
\smallskip

\noindent{\it Proof of Proposition \ref{dif} for $b_D.$}
Let $\eps\in(0,1/4].$ The arguments from the previous proof and Lemma L implies
that there exist neighborhoods $U_2\subset U_1$ of $a$ such that
$F=D\cap U_1$ is a Dini-smooth simply connected domain and
$\beta_D>\beta_F-\eps$ on $D\cap U_2.$

Let $z,w\in D.$ We may choose a smooth curve
$\gamma:[0,1]\to D$ with $\gamma(0)=z,$ $\gamma(1)=w$ and
$$b_D(z,w)+\eps>\int_0^1(\beta_D(\gamma(t))|d\gamma(t)|.$$
It follows by the proof of \cite[Proposition 1]{NT} that one may find
a neighborhood $U_3\subset U_2$ of $a$ such that $\gamma\subset D\cap U_2$ if
$z,w\in D\cap U_3.$

Let $\psi:F\to\D$ be a Riemann map, $\tilde\gamma=\psi\circ\gamma,$ $z,w\in D\cap U_3$
and $\tilde z=\psi(z),$ $\tilde w=\psi(w).$ Using that $\eps\le 1/4,$ we obtain that
\begin{align*}\sqrt2(b_D(z,w)+\eps)&>\sqrt2 \int_0^1(\beta_\D(\psi(t))-\eps)|d\psi(t)|>\int_0^1\frac{|d\psi(t)|}
{d_\D(\psi(t))}\\
&\ge h_\D(\tilde z,\tilde w)\ge 2k_\D(\tilde z,\tilde w)-c|\tilde z-\tilde w|=:2\tilde s_D(z,w),
\end{align*}
where $c>0$ is a constant. Since $\ds\lim_{z,w\to a}(s_D(z,w)-\tilde s_D(z,w))=0,$ we get
$$\liminf_{z,w\to a}(b_D(z,w)-\sqrt2s_D(z,w))\ge 0.$$

The proof of the opposite inequality
\begin{equation}\label{opp}
\limsup_{z,w\to a}(b_D(z,w)-\sqrt2s_D(z,w))\le 0.
\end{equation}
is even simpler. We choose $U_2$ such that $\beta_D<\beta_F+\eps$ on $D\cap U_2.$
Then we may take $U_3$ such that the $b_F$-geodesic for any $z,w\in D\cap U_3$ belongs to
$D\cap U_2.$ It follows that $b_D(z,w)\le b_F(z,w)+c'|z-w|$ for some constant $c'>0$ which implies \eqref{opp}.
\smallskip

\noindent{\it Proof of Proposition \ref{rat}.} (a) First, we shall show that
$$\lim_{\substack{z,w\to a\\z\neq w}}\frac{c_D(z,w)}{k_D(z,w)}=1.$$

We may proceed as in the proof of Proposition \ref{dif}, choosing $\zeta$
to be a $C^1$-smooth curve and then observing that $\psi$ extends to a homeomorphism
from $\overline G$ to $\overline \D.$ Then there exists a disc $U$ centered at $\eta(a)$
such that $\eta(D)\subset\D\cap U.$ Hence
$$1\ge\frac{c_D(z,w)}{k_D(z,w)}\ge\frac{c_\D(\eta(z),\eta(w))}
{k_{\D\cap U}(\eta(z),\eta(w))},\quad z\neq w.$$
It follows similarly to the proof of Lemma L that the last quotient tends to 1 as $z,w\to a.$

Second, we have that (see \cite{JN})
$$\lim_{z\to a}2\kappa_D(z)d_D(z)=\sqrt2\beta_D(z)d_D(z)=1.$$
It follows as in the proof of \cite[Proposition 5]{NA} that
$$\lim_{\substack{z,w\to a\\z\neq w}}\frac{h_D(z,w)}{p_D(z,w)}=2.$$
It remains to use that, by \cite[Proposition 6 (a)]{NA},
$$\lim_{\substack{z,w\to a\\z\neq w}}\frac{h_D(z,w)}{s_D(z,w)}=2.$$

\noindent(b) It is a direct consequence of (a) for $k_D$ and $b_D,$ since they are inner distances.
On the other hand, by \cite[Proposition 9]{N1},
$$\lim_{z\to\partial D}\frac{c_D(z,w)}{k_D(z,w)}=1
\quad\mbox{ uniformly in }w\neq z$$
which completes the proof.

\end{document}